\documentclass[11pt,a4paper,english]{amsart}
\usepackage[english]{babel}
\usepackage[latin1]{inputenc}
\usepackage{amssymb}
\usepackage{amsmath}
\usepackage{amsthm}
\usepackage{bbm}

\newcommand\RR{\mathbb{R}}
\newcommand\Sph{\mathbb{S}}
\newcommand\OO{\ensuremath{\Omega}}
\newcommand\CC{\mathbb{C}}

\newcommand\tr{\ensuremath{\triangle}}
\newcommand\Str{Strichartz }
\newcommand\Schr{Schrödinger }
\newcommand\GP{Gross-Pitaevskii }

\newcommand\pH{\ensuremath{\dot{H}^1(\Omega)}}

\newcommand\norm[2]{\ensuremath{|\!| #1 |\!|_{#2}}}

\newcommand\nL[2]{\norm{ #1}{L^{#2}}}
\newcommand\nLi[3]{\norm{#1}{L^#2(#3)}}
\newcommand\nLdi[2]{\norm{#1}{L^2(#2)}}
\newcommand\nLdO[1]{\norm{#1}{L^2(\Omega)}}

\newcommand\nH[2]{\norm{#1}{H^{#2}}}
\newcommand\nHO[2]{\norm{#1}{H^{#2}(\Omega)}}

\newcommand\nX[2]{\norm{#1}{#2}}
\newcommand\nLL[3]{\norm{#1}{L^{#2}(L^{#3})}}

\newtheorem{Def}{Definition}
\newtheorem{thm}{Theorem}[section]
\newtheorem{lem}[thm]{Lemma}
\newtheorem{propos}[thm]{Proposition}
\newtheorem{rem}{Remark}
\newtheorem{cor}[thm]{Corollary}

\newtheorem*{citelem}{Lemma}
\newtheorem*{citeprop}{Proposition}

\title[NLS and GP in exterior domains]{Global existence for defocusing cubic NLS and Gross-Pitaevskii equations in three dimensional exterior domains}
\author{Ramona Anton}
\address{Universit\'e Paris Sud, Math\'ematiques, B\^at 425, 91405 Orsay Cedex}
\email{Ramona.anton@math.u-psud.fr}
\date{}
\begin{document}

\begin{abstract}We prove global wellposedness in the energy space of the defocusing cubic nonlinear Schrödinger and Gross-Pitaevskii equations on the exterior of a non-trapping domain in dimension 3. The main ingredient is a Strichartz estimate obtained combining a semi-classical Strichartz estimate \cite{RA} with a smoothing effect on exterior domains \cite{BGTannIHP}.

\end{abstract}

\maketitle

\section{Introduction}

Let $\Theta\neq \emptyset$, $\Theta\subset \RR^3$, a nontrapping obstacle with compact boundary and let $\OO=\complement \Theta$. In this paper we are interested in the Cauchy problem for the cubic defocusing NLS equation (here written with Dirichlet boundary conditions) on $\OO$ : 
\begin{equation}
\label{CubicNLS}
\left \{
\begin{array}{rcl}
i \partial_t u + \triangle u &=& |u|^2 u, \ on \ \mathbb{R} \times
\Omega \\ u_{|_{t=0}} &=& u_0, \ on \ \Omega\\
u_{|_{\mathbb R \times \partial \Omega}} &=& 0.
\end{array}
\right.
\end{equation}
This equation appears in the nonlinear optics and more generally in propagation of nonlinear waves. For more details on nonlinear \Schr equations see for example the  books of C.Sulem-P.L.Sulem \cite{SuSu}, T.Cazenave \cite{Ca} and the references therein.

There is a wide literature on the Cauchy problem in the Euclidean space. One of the main tools in addressing this problem is the Strichartz inequality, which translates the dispersive property of the linear \Schr flow. We refer to the work of Strichartz \cite{STR}, Ginibre-Velo \cite{GiVe} and Keel-Tao \cite{KT}. 

Recently, the question of the influence of the geometry on the solution has been studied. Let us mention the work of J.Bourgain \cite{Bo} on the tori $\mathbb{T}^d$ for $d=2,3$ and of N.Burq-P.Gérard-N.Tzvetkov \cite{BGTajm}, \cite{BGTannIHP} on compact manifold and exterior of non-trapping obstacles.

In recent works on superfluidity and Bose-Einstein condensates (see for example the book of A.Aftalion \cite{Afta}) the following variant of NLS (\ref{CubicNLS}) is studied 
\begin{equation}
\label{CubicGP}
\left \{
\begin{array}{rcl}
i \partial_t u + \triangle u &=& (|u|^2-1) u, \ on \ \mathbb{R} \times
\Omega \\ u_{|_{t=0}} &=& u_0, \ on \ \Omega\\
\frac{\partial u}{\partial \nu}{|_{\mathbb R \times \partial \Omega}} &=& 0.
\end{array}
\right.
\end{equation}
This is called the cubic Gross-Pitaevskii equation with Neumann boundary conditions. 
The main difference between the NLS (\ref{CubicNLS}) and the \GP equation (\ref{CubicGP}) is in their energy space. For \GP it reads
$$E=\{u\in H_{loc}^1(\Omega),\ \nabla u \in L^2(\Omega),\ |u|^2-1 \in L^2(\Omega)\}.$$
Namely, the initial datum in the energy space, $u_0\in E$, is not an $L^2(\Omega)$ function.
In \cite{BeSa}, \cite{BeOS}, \cite{AftaBl}, \cite{Gravj1}, \cite{Gravj2}, \cite{Chiron} the question of existence of travelling waves and vortices is studied. We are interested here in providing a mathematical background for the study of the dynamics of these phenomena. More precisely, we are interested in showing wellposedness in the energy space. There have been previous works on the Cauchy problem for the \GP equation : P.E.Zhidkov \cite{Zdkv1}, \cite{Zdkv2} in Zhidkov spaces $X^1(\RR)$, F.Béthuel- J.C.Saut \cite{BeSa} in the space of functions $1+H^1(\RR^d)$, for $d=2,3$, P.Gérard in \cite{PG} in the energy space on the whole Euclidean space $\RR^d$, for $d=2,3,4$, C.Gallo \cite{Ga} in the energy space $u_0+H^1(\OO)$ for exterior domains in $d=2$.

For both (\ref{CubicNLS}) and (\ref{CubicGP}) the method we use is based on a new Strichartz estimate obtained combining a smoothing effect in exterior domains \cite{BGTannIHP} with a semiclassical Strichartz estimate on small intervals of time depending on the frequencies where the flow is localised \cite{RA}. In dimension $2$ the smoothing effect \cite{BGTannIHP} provides wellposedness for both NLS \cite{BGTannIHP} and \GP \cite{Ga}, with all power nonlinearities. In dimension $3$ the smoothing effect only provides wellposedness of subcubic nonlinearities \cite{BGTannIHP}, \cite{Ga}. Improving the Strichartz inequality allows us to treat the cubic nonlinearity in exterior domains in dimension $3$.

Let us recall the definition of an admissible pair. 

\begin{Def} A pair $(p,q)$ is called admissible in dimension $3$ if $p\geq 2$ and $$\frac 2p +\frac 3q=\frac 32.$$
\end{Def}

The \Str inequality we obtain is the following.

\begin{propos}
\label{StrEstExtDom}
For $(p,q)$ an admissible couple in dimension $3$ and $\epsilon>0$, there exists $c_\epsilon>0$ such that for all $u_0\in H_D^{\frac {1}{2p} + \epsilon}(\Omega)$,
\begin{equation}
\label{StrExt}
\nX{e^{it\tr_D} u_0}{L^p(I,L^q(\Omega))} \leq c_\epsilon \nX{u_0}{H_D^{\frac {1}{2p} + \epsilon}(\Omega)}
\end{equation}
A similar result holds for the linear \Schr flow with Neumann boundary conditions.
\end{propos}

Having a \Str inequality we obtain classically a local existence theorem for (\ref{CubicNLS}) by Picard iteration scheme. These also enables propagation of the regularity of the initial data. Local existence in the energy space $H_0^1(\OO)$ combined with the conservation of the energy (and for defocusing nonlinearity of the $H_0^1(\OO)$ norm) enables us to conclude that the solution to (\ref{CubicNLS}) is global in time.

\begin{thm}
\label{LeNLS}
For all $u_0\in H_0^1(\Omega)$ there exists an unique solution $$u \in C(\RR, H_0^1(\Omega)) \cap L_{loc}^p(\RR,L^\infty(\Omega))$$ (for every $2<p<3$) of equation (\ref{CubicNLS}). Moreover, for every $T>0$ and for every bounded subset $B$ of $H_0^1(\Omega)$, the flow $u_0 \mapsto u$ is Lipschitz from $B$ to $C([-T,T],H_0^1(\Omega))$. For $1<\sigma\leq 2$ and $u_0\in H_0^1(\OO) \cap H^\sigma (\OO)$ we have $u\in C([-T,T], H_0^1(\OO) \cap H^\sigma (\OO))$. 
\end{thm}

For \GP equation, as $u_0\in E$ is not an $L^2(\OO)$ function, the \Str inequality does not apply directly. We adapt the arguments of \cite{PG} to the boundary case for the description of the structure of $E$ and of the action of the linear \Schr group on $E$. In particular we define a structure of complete metric space on $E$. The global existence theorem for the \GP equation (\ref{CubicGP}) follows by combining the latter structure with dispersive estimates (\ref{StrExt}).
\begin{thm}
\label{GeGP} 
For all $u_0\in E$ there exists an unique solution $$u \in C(\RR, E) \cap L_{loc}^p(\RR,L^\infty(\Omega))$$ (for every $2<p<3$) of equation (\ref{CubicGP}). Moreover, the following properties hold: for every bounded subset $B$ of $E$ there exists $T>0$ such that for all $u_0 \in B$ the flow $u_0 \mapsto u$ is Lipschitz from $B$ to $C([-T,T],E)$ ; we have $u-u_L\in C(\RR,H^1(\OO))$ ; if $u_0\in E$ is such that $\tr u_0 \in L^2(\OO)$ and $\frac{\partial u_0}{\partial \nu}=0$, then $\tr u \in C(\RR, L^2(\OO))$. 
\end{thm}

\begin{rem}
After the completion of this work Blair-Smith-Sogge \cite{BlSmSo} announced an improved \Str inequality on boundary domains. They prove a \Str inequality with a loss of $\frac{4}{3p}$ derivatives as opposed to the \Str inequality \cite{RA} with a loss of $\frac{3}{2p}$ derivatives we used here. Although this may improve our \Str inequality (\ref{StrExt}) it does not improve the wellposedness results of Theorem \ref{LeNLS} and of Theorem \ref{GeGP}. 
\end{rem}

The structure of the paper is as follows : in Section \ref{StrSect} we show how we obtain the \Str estimate (\ref{StrExt}). In Section \ref{NLSsect} we give the proof of Theorem \ref{LeNLS}. In Section {\ref{GPsect}} we deal with the \GP equation (\ref{CubicNLS}) and we give the proof of Theorem \ref{GeGP}.

\noindent {\bf Acknowledgments :} \textit{The author would like to thank P.Gérard for proposing the subject of this article and for many helpfull and stimulating discussions. The result of this article form part of author's PhD thesis defended at Université Paris Sud, Orsay, under P.Gérard's direction.}


\section{Strichartz estimate in exterior domains}
\label{StrSect}
The idea is to combine \Str inequality on exterior domains \cite{RA} 
with the gain of $\frac 12$ derivative from the smoothing effect \cite{BGTannIHP}. Rather than using the Strichartz estimate with loss of $\frac{3}{2p}+\epsilon$ derivatives (45) of \cite{RA}, we prefer to use the Strichartz estimate without loss of derivatives (Proposition 4.13 of \cite{RA}), that holds for frequency localised initial data and small intervals of time depending on the frequency.

In order to do that here, we need to recall some of the notations and results from \cite{RA}. That is done in Subsection \ref{Pre}. In Subsection \ref{rappel} we recall the results of N.Burq, P.Gérard and N.Tzvetkov \cite{BGTannIHP} concerning the smoothing effect and \Str estimate away from the obstacle. Subsection \ref{near} is the core of this section. We prove a new \Str estimate close to the obstacle by combining semiclassical \Str estimate and smoothing effect. In Subsection \ref{StrPr} we deduce the proof of Proposition \ref{StrEstExtDom}.

\subsection{Preliminaries}
\label{Pre}

We recall here the classical mirror reflection that allows us to pass from a manifold with boundary to a boundaryless manifold. This method consists in taking a copy of the domain and glue it to the initial one by identifying the points of the boundary. In order for this to be a manifold we have to choose the coordinates carefully. Thus, taking normal coordinates at the boundary is like straightening a neighborhood of the boundary into a cylinder $\partial\Omega \times [0,1)$ and gluing the two cylinders along the boundary makes a nice smooth manifold. This can be properly done using for example tubular neighborhoods (e.g. \cite{Spiv}, pp. 468 and 74). Let $M=\Omega\times\{0\}\cup_{\partial\Omega}\Omega\times\{1\}$, where we identify $(p,0)$ with $(p,1)$ for $p\in \partial\Omega$. 
\begin{citelem} \textbf{\upshape{(\cite{Spiv})}}There is a unique $C^\infty$ structure on M such that $\Omega\times \{j\} \hookrightarrow M$ is $C^\infty$ and $\tilde{\chi}: U\times\{0\} \cup_{\partial\Omega} U\times\{1\} \rightarrow \partial\Omega \times (-1,1)$ is a diffeomorphism, where $U$ is a small neighborhood of $\partial \Omega$ for which there are deformation retractions onto $\partial \Omega$.
\end{citelem}
On $M$ we define the metric $G$ induced by the new coordinates. As we have chosen coordinates in the normal direction, the metric is well defined over the boundary, its coefficients are Lipschitz in local coordinates and diagonal by blocs (no interaction between the normal and the tangent components). Moreover, $$G(r(y))=G(y),$$ where $r:M\rightarrow M,\ r(x,0)=(x,1),\ r^2=Id$ is the reflection with respect to the boundary $\partial\Omega$.

For the {\bf Dirichlet} problem we introduce the space $H_{AS}^1$ of functions of $H^1(M)$ which are anti-symmetric with respect to the boundary. Let
$$H_{AS}^1=\{v:M\rightarrow \mathbb{C},\ v\in H^1(M),\ v(y)=-v(r(y))\}.$$

\noindent Note that for $v\in H_{AS}^1$ the restriction $v_{|_{\Omega\times{0}}}$ is in $H_0^1(\Omega)$ and every function from $H_{AS}^1$ is obtained from a function of $H_0^1(\Omega)$. We shall prove the stability of $H_{AS}^1$ under the action of $e^{it\tr_G}$. 

By complex interpolation define $H_{AS}^s$ for $s\in [0,1]$ and deduce its stability under the action of $e^{it\tr_G}$. Moreover, the restriction to $\Omega$ of functions in $H_{AS}^s$ belongs to $H_D^s(\Omega)$ and vice versa. This allows us to deduce the Strichartz inequality for $e^{it\tr_D}$ on $\Omega$ from the Strichartz inequality for $e^{it\tr_G}$ on $M$. 

Similarly, we can define for the {\bf Neumann} problem the space $H_S^1$ of symmetric functions with respect to the boundary. This space is also stable under the action of $e^{it\tr_G}$.
$$H_{S}^1=\{v:M\rightarrow \mathbb{C},\ v \in H^1(M),\ v(y)=v(r(y))\}.$$

\label{reflex}
Let us prove the stability of $H_{AS}^1$ under the action of $e^{it\tr_G}$. Let $v_0\in H_{AS}^1$ and $v(t,y)=e^{it\tr_G} v_0$. Then $v$ satisfies to $i \partial_t v(t,y)+ \tr_{G(y)} v (t,y)= 0$, $v(0)=v_0$. Let $\tilde{v}(t,y)=v(t,r(y)).$ We shall look for the equation verified by $\tilde{v}$. First note that $\tilde{v}(0)=-v_0$ and $\partial_t \tilde{v}(t,y)=\partial_t v(t,y)$. As $G$ is diagonal by blocks, having no interactions between the normal and tangent components, so is $G^{-1}$. Thus in $\tr_{G(y)}$ there is no crossed term. Consequently $\tr_{G(r(y))}\tilde{v}(t,y) = \tr_{G(y)}v(t,y)$. We see thus that $\tilde{v}$ satisfies to the linear Schrödinger equation with initial data $-v_0(y)$. But $-v(t,y)$ satisfies the same equations. By uniqueness we conclude that $$v(t,r(y))=-v(t,y).$$
Moreover, if $v_0(t,y)=u_0(t,y)$ for all $y\in \Omega$, then $v(t,y)=u(t,y)$ for all $t$ and for all $y\in \Omega$, where $u(t)=e^{it\tr_D}u_0$.

We prepare the frequency decomposition. We begin with a partition of unity on $M$. Since $M$ is flat outside a compact set, let $(U_j,\kappa_j)_{j\in J}$ be a covering of the area of $M$ where $G\neq \mathbb{I}d$. This area is compact, so we can choose $J$ of finite cardinal. We have $M=\cup_{j\in J}U_j \cup U_{1,\infty}\cup U_{2,\infty}$, where $U_{1,\infty}$ and $U_{2,\infty}$ are two disjoint neighborhood of $\infty$, diffeomorphe to $\RR^d \backslash \bar{B}$. Let $(\chi_j)_{j\in J},\ \chi_{1,\infty},\ \chi_{2,\infty} : M \rightarrow [0,1]$ be a partition of unity subordinated to the previous covering. For all $j \in J$ let $\tilde{\chi}_j : M \rightarrow [0,1]$ be a $C^\infty$ function such that $\tilde\chi_j=1$ on the support of $\chi_j$ and the support of $\tilde\chi_j$ is contained in $U_j$. Similarly we define $\tilde{\chi}_{1,\infty},\ \tilde{\chi}_{2,\infty} : M \rightarrow [0,1]$. Let $\varphi_0 \in C^\infty(\RR^d)$ be supported in a ball centered at origin and $\varphi \in C^\infty(\RR^d)$ be supported in an annulus such that for all $\lambda \in \RR^d$
\begin{equation}
\label{LPdecomp}
\varphi_0(\lambda) + \sum_{k \in \mathbb{N}} \varphi(2^{-k} \lambda)=1.
\end{equation}

We define a family of spectral truncations : for $f\in C^\infty(M)$ and $h\in (0,1)$ let
\begin{equation}
\label{defJh}
J_h f=\sum_{j\in J} (\kappa_j)^* \left( \tilde{\chi}_j \varphi(hD)(\kappa_j^{-1})^* (\chi_j f) \right) + F_{1,h,\infty}f + F_{2,h,\infty}f
\end{equation} and
\begin{equation}
J_0 f=\sum_{j\in J} (\kappa_j)^* \left( \tilde{\chi}_j \varphi_0(D)(\kappa_j^{-1})^* (\chi_j f) \right) +F_{1,0,\infty}f+F_{2,0,\infty}f,
\end{equation}
where $^*$ denotes the usual pullback operation and $F_{h,\infty} f =\tilde{\chi}_{\infty} \varphi(hD) \chi_{\infty}f(x).$ The following identity holds
$$ J_0 f(x) + \sum_{k=0}^\infty J_{2^{-k}} f(x) = f(x).$$
This will be useful for the Littlewood Paley theory.

We need more regularity on the coefficients of the metric than the Lipschitz regularity. Therefore we define a regularized metric $G_h$ as follows : let $\psi$ be a $C_0^\infty(\RR^d)$ radially symmetric function with $\psi\equiv 1$ near 0. Let
\begin{equation}
\label{regmetricdef}
G_h =\sum_{j\in J} (\kappa_j)^* \left( \tilde{\chi}_j \psi(h^{\frac 12} D)(\kappa_j^{-1})^* (\chi_j G) \right).
\end{equation}
The transformation of $G$ into $G_h$ does not spoil the symmetry. Note also that $G_h$ converges uniformly in $x$ to $G$, and thus, for $h$ sufficiently small, $G_h$ is positive definite. Therefore, $G_h$ is still a metric. We present some properties of metrics $G$ and $G_h$.

\begin{citelem} The metric $G:M\rightarrow M_d(\RR)$ is symmetric, positive definite and Lipschitz : there exist $c,C,c_1 >0$ such that for all $x\in M$
$$c{\mathbb{I} d} \leq G(x) \leq C {\mathbb{I} d},\ |\partial G|\leq c_1,$$
where we have denoted by $\partial G$ the derivatives of the metric in a system of coordinates. The regularized metric $G_h$ is a $C^\infty$ function
that verifies the followings : there exists $c, C>0$ and $c_\gamma >0$ for all
$\gamma \in \mathbb{N}^d$ such that 
$$c {\mathbb I}d \leq G_h(x) \leq C {\mathbb I}d,\ |\partial^\gamma G_h(x)| \leq c_\gamma h^{-\alpha \max (|\gamma|-1 , 0)}.$$
\end{citelem}

We present next a collection of estimates on $J_h$. There exist constants $c>0$ such that, for all $h\in (0,1)$ :
\begin{itemize}
	\item $\nX{J_h}{L^p \rightarrow L^p} \leq c_p$, for all $1\leq p \leq \infty$.
	\item $\nX{[J_h,\tr_{G_h}]}{L^2\rightarrow L^2} \leq \frac{c}{h}$ and  $\nX{[J_h,\tr_{G_h}]}{H^1\rightarrow L^2} \leq c$. 
	
	As one may not apply two derivatives on $G(x)$, the similar statement for $[F_h,\tr_G]$ only holds for the $H^1 \rightarrow L^2$ norm : $\nX{[F_h,\tr_G]}{H^1\rightarrow L^2} \leq c.$
	\item $\nX{J_h(\tr_{G_h}-\tr_G)}{H^1 \rightarrow L^2} \leq c h^{-\frac 12}.$
\end{itemize}

We define also a spectral cut-off slightly larger than $J_h$. 
Let $\tilde{\varphi}$ be a $C^\infty$ function supported in an annulus such that $\tilde{\varphi}=1$ on a neighborhood of the support of $\varphi$. We define $\tilde{J}_h$ just like $J_h$, replacing $\varphi$ par $\tilde{\varphi}$ in (\ref{defJh}) :
\begin{equation}
\label{art2defTildeJh}
\tilde{J}_h f=\sum_{j\in J} (\kappa_j)^* \left( \tilde{\chi}_j \tilde \varphi(hD)(\kappa_j^{-1})^* (\chi_j f) \right) + \tilde F_{1,h,\infty}f + \tilde F_{2,h,\infty}f
\end{equation}
Then the action of $\tilde{J}_h$ on $J_h$ and $[J_h, \tr_{G_h}]$ is close to identity in $L^p\rightarrow L^p$ norm, $p\geq 2$, and $L^2 \rightarrow L^2$ norm respectively.
\begin{itemize}
  \item $\nX{\tilde{J}_h J_h - J_h}{L^p \rightarrow L^p} \leq c_N h^N$,
	\item $\nX{[J_h,\tr_{G_h}]- [J_h,\tr_{G_h}]\tilde{J}_h}{L^2\rightarrow L^2} \leq c_N h^N$, for all $N\in \mathbb{N}$.
\end{itemize}

Let us recall also the Strichartz estimate we use from \cite{RA}.

\begin{citelem} (4.13 of \cite{RA}) For all couples $(p,q)$ admissible in dimension $3$ and $I_h$ an interval of time such that $|I_h|=ch^{\frac 32}$, we have
\begin{equation}
\label{hStrinter}
\nX{J_h^* e^{it\tr_{G_h}} u_0}{L^p(I_h,L^q(M))} \leq c\nL{u_0}{2}
\end{equation}
\end{citelem}

We prefer to go back to the estimate on $e^{it\tr_{G_h}}$ since the form of the Strichartz estimate for $e^{it\tr_G}$ is more difficult to handle ((45) of \cite{RA}) :
$$\nX{J_h^* e^{it\tr_{G}} u_0}{L^p(I_h,L^q(M))} \leq ch \nH{u_0}{1}.$$

This is due to the fact that $\tr_G$ and $\tr_{G_h}$ are not both selfadjoint in the same space, because of the  the volume density $\frac{1} {\sqrt{\det{G(x)}}}$.

\subsection{Smoothing effect and \Str estimate away from the obstacle}
\label{rappel}
In this section we recall two results of N.Burq, P.Gérard and N.Tzvetkov \cite{BGTannIHP} on the smoothing effect for the \Schr flow on exterior domains and the Strichartz estimate away from the obstacle. The smoothing effect was obtained via resolvent bounds. For \Str estimate they used a strategy inspired by G.Staffilani and D.Tataru's paper \cite{StTa} on $C^2$ short range perturbation of the free Laplacian on $\RR^d$. Thus, they proved that away from the obstacle the linear \Schr flow satisfies the usual \Str estimates. We present an equivalent statement on the double manifold. 

\begin{citeprop} (2.7 of \cite{BGTannIHP})
Assume that $\Theta\neq \emptyset$. Then for every $T>0$, for every $\chi\in C_0^\infty(\RR^d)$, $d\geq 2$,
$$
\nX{\chi e^{i t\tr_D} u_0}{L_T^2H_D^{s+ \frac 12}(\Omega)} \leq c \nX{u_0}{L_T^2H_D^{s}(\Omega)}, \ {\rm for} \ s\in [0,1]. 
$$
\end{citeprop}

\begin{citeprop} (2.10 of \cite{BGTannIHP}) For every $T>0$, for every $\chi\in C_0^\infty(\RR^d)$, $\chi=1$ close to $\Theta$, there exists $C>0$ such that
\begin{equation}
\label{StrCloseInit}
\nX{(1-\chi)u}{L_T^pW^{s,q}(\Omega)} \leq C\nX{u_0}{H_D^s}(\Omega),
\end{equation}
where $s\in [0,1]$, $u(t)=e^{it\tr_D}u_0$ and $(p,q)$ any \Str admissible pair.
\end{citeprop}

%

The proof relies on the use of the smoothing effect and the fact that $(1-\chi)e^{it\tr_D}u_0$ can be seen as a solution to some nonlinear \Schr equation on $\RR^d$.

Although the properties are written for the Dirichlet Laplacian, Remark 1.2 of \cite{BGTannIHP} ensures that the results hold for the Neumann conditions as well. From the way we constructed the double manifold and flows, we deduce that those results extend easily on the double manifold.

\begin{propos}
Assume that $\Theta\neq \emptyset$. Then for every $T>0$, for every $\chi\in C_0^\infty(M)$, 
$$
\nX{\chi e^{i t\tr_G} u_0}{L_T^2H^{s+ \frac 12}(M)} \leq c \nX{u_0}{L_T^2H^{s}(M)}, \ {\rm for} \ s\in [0,1]. 
$$
\end{propos}

\begin{propos}
For every $T>0$, for every $\chi\in C_0^\infty(M)$, $\chi=1$ close to ${\mathcal D}\Theta$, where ${\mathcal D}\Theta$ represents the double of $\Theta$, there exists $C>0$ such that
\begin{equation}
\label{StrAway}
\nX{(1-\chi)e^{it\tr_G}u_0}{L_T^p W^{s,q}(M)} \leq C\nX{u_0}{H^s(M)},
\end{equation}
where $s\in [0,1]$ and $(p,q)$ any \Str admissible pair.
\end{propos}

\subsection{\Str estimate near the obstacle}
\label{near}
We want to combine Strichartz estimate on domains \cite{RA} with smoothing effect \cite{BGTannIHP}. For this we use the Strichartz estimate of the frequency localised linear flow, without loss of derivatives, which holds on a small interval of time (see estimate (\ref{hStrinter}) from Subsection \ref{Pre}).

Let $\varphi \in C_0^\infty(\RR)$ such that $\varphi\equiv 1$ on $[-\frac 14, \frac 14]$, $0 \leq \varphi \leq 1$, $supp \varphi \subset [-\frac 12, \frac 12]$ and there exists $J\subset \RR$ a discrete set such that $\sum_{t_0 \in J} \varphi^2(t-t_0) = 1$ for all $t\in \RR$. Let $\delta>0$ be a small number. If we consider $\bar{J}=[-\delta, 1+\delta] \cap ch^{\frac 32} J$ then,  for $t\in [-\frac \delta 2, 1+\frac \delta 2]$, 
\begin{equation}
\label{idPhi01}
\sum_{t_0 \in \bar{J}} \varphi^2 \left( \frac{t-t_0}{ch^{\frac 32}} \right) = 1.
\end{equation}
Let us denote by $I_h(t_0)=[t_0-\frac{ch^{\frac 32}}{4}, t_0+\frac{ch^{\frac 32}}{4}]$, $I'_h(t_0)=[t_0-\frac{ch^{\frac 32}}{2}, t_0+\frac{ch^{\frac 32}}{2}]$,  $u_L(t)=e^{it\tr_G}u_0$ and by 
\begin{equation}
\label{defv}
v(t)=\varphi\left( \frac{t-t_0}{ch^{\frac 32}} \right) J_h^* \chi e^{it\tr_G} u_0.
\end{equation}
Notice that $v(t)=J_h^* \chi e^{it\tr_G} u_0$ for $t\in I_h(t_0)$ and $supp_t v \subset I'_h$. We write the end-point Strichartz estimate for $v$ on $I'_h$. Notice that the couple $(2,6)$ is admissible in dimension $3$.

\begin{lem} 
\label{1truncStr}
For $t_0\in \RR$, $\tilde{\chi}\in C_0^\infty(M)$ such that $\tilde{\chi}\chi = \chi$ and $\tilde{J}_h$ a spectral cut-off slightly larger than $J_h$ (\ref{art2defTildeJh}), 
we have 
\begin{eqnarray}
\label{trunc1}
\nLdi{\varphi\left( \frac{t-t_0}{ch^{\frac 32}} \right) J_h^* \chi e^{it\tr_G} u_0}{I'_h, L^6(M)} \leq  ch^{-\frac 34} \nLdi{J_h^* \chi u_L}{I'_h, L^2(M)} + \nonumber \\h^{\frac 34} \nLdi{\tilde{J}_h^* \tilde{\chi} u_L}{I'_h, H^1(M)} + c\nLdi{\tilde{\chi} u_L}{I'_h, H^1(M)}.
\end{eqnarray}
\end{lem}

\begin{proof}
For simplicity, let us suppose that $I_h'(t_0)=[0,T]$, where $T=ch^{\frac 32}$. Then $v(t)$ verifies, for $t\in I_h'(t_0)$, the equation
$$
\left \{
\begin{array}{rcl}
i \partial_t v + \triangle_{G_h} v &=& f_1 + f_2 +f_3\\
v_{|_{t=0}} &=& 0
\end{array}
\right.
$$
where $f_1=\frac{i}{ch^{\frac 32}} \varphi'\left( \frac{t-t_0}{ch^{\frac 32}} \right) J_h^* \chi u_L$, $f_2= \varphi\left( \frac{t-t_0}{ch^{\frac 32}} \right) [\tr_{G_h}, J_h^* \chi] u_L $ and $f_3=\varphi\left( \frac{t-t_0}{ch^{\frac 32}} \right) J_h^* \chi (\tr_G - \tr_{G_h}) u_L$. By the Duhamel formula and using that $\tilde{J}_h^* J_h^* = J_h^* + c_N h^N$ in $L^p\rightarrow L^p$ norm, for $p\geq 2$, we have $$v(t)=v_1(t)+ v_2(t)+ v_3(t) + c_N h^N,$$ where 
we define
$$
\begin{array}{rcl}
v_1(t) &=& \frac{1}{ch^{\frac 32}} \int_0^t \tilde{J}_h^* e^{i(t-\tau)\tr_{G_h}} \varphi'\left( \frac{\tau -t_0} {ch^{\frac 32}} \right) J_h^* \chi u_L(\tau) d\tau,\\
v_2(t) &=& -i \int_0^t \tilde{J}_h^* e^{i(t-\tau)\tr_{G_h}} \varphi\left( \frac{\tau -t_0}{ch^{\frac 32}} \right) [\tr_{G_h}, J_h^* \chi] u_L(\tau) d\tau,\\
v_3(t) &=& -i \int_0^t \tilde{J}_h^* e^{i(t-\tau)\tr_{G_h}} \varphi\left( \frac{\tau - t_0}{ch^{\frac 32}} \right) J_h^* \chi (\tr_G -\tr_{G_h}) u_L(\tau) d\tau.
\end{array}
$$

By Minkowski inequality and estimate (\ref{hStrinter}), we have 
$$
\begin{array}{rcl}
\nX{v_1}{L^2_t(I_h', L^6(M))} &\leq & ch^{-\frac 32} \int_0^T |\varphi'\left( \frac{\tau -t_0} {ch^{\frac 32}} \right) | \nX{\mathbbm{1}_{\tau<t} \tilde{J}_h^* e^{i(t-\tau)\tr_{G_h}} J_h^* \chi u_L(\tau)}{L^2_t(L_x^6(M))} d\tau, \\
& \leq & ch^{-\frac 32} \int_0^T |\varphi'\left( \frac{\tau -t_0} {ch^{\frac 32}} \right) | \nX{J_h^* \chi u_L(\tau)}{L^2_x(M)} d\tau.
\end{array}
$$
Using Cauchy Schwarz inequality and $\nL{\varphi'(\frac{\cdot}{ch^\frac 32})}{2} =ch^{\frac 34}$, we obtain
\begin{equation}
\label{StrV1}
\nX{v_1}{L^2_t(I_h', L^6(M))} \leq ch^{-\frac 34} \nX{J_h^* \chi u_L}{L^2(I_h'\times M)}
\end{equation}

\noindent Similarly, we have $\nX{v_2}{L^2_t(I_h', L^6(M))} \leq ch^{ \frac 34} \nX{[\tr_{G_h},J_h^* \chi] u_L}{L^2(I_h'\times M)}.$
Using that $\nX{[\tr_{G_h},J_h^* \chi]}{H^1 \rightarrow L^2} \leq c$ and $\nX{[\tr_{G_h},J_h^* \chi]}{H^1 \rightarrow L^2} \sim \nX{[\tr_{G_h},J_h^* \chi]\tilde{J}_h^* \tilde{\chi}}{H^1 \rightarrow L^2}$ modulo $ch^{-N}$, we obtain
\begin{equation}
\label{StrV2}
\nX{v_2}{L^2_t(I_h', L^6(M))} \leq ch^{\frac 34} \nX{\tilde{J}_h^* \tilde{\chi} u_L}{L^2(I_h', H^1(M))}.
\end{equation}

We estimate the third term $v_3$ in $L^2(I'_h,L^6(M))$ norm in a similar manner. We get : $\nX{v_3}{L^2_t(I_h', L^6(M))} \leq ch^{ \frac 34} \nX{J_h^* \chi (\tr_G - \tr_{G_h}) u_L}{L^2(I_h'\times M)}.$ Using the estimate $\nLdi{J_h^* \chi (\tr_G - \tr_{G_h}) f}{M} \leq ch^{-\frac 12} \nX{\tilde{\chi} f}{H^1(M)}$, we obtain
\begin{equation}
\label{StrV3}
\nX{v_3}{L^2_t(I_h', L^6(M))} \leq ch^{\frac 14} \nX{\tilde{\chi} u_L}{L^2(I_h', H^1(M))}.
\end{equation}

Recalling that $v(t)=v_1(t)+ v_2(t)+ v_3(t)$, the result follows from the triangle inequality and the sum of (\ref{StrV1}), (\ref{StrV2}) and (\ref{StrV3}).
\end{proof}

We proceed to the summation over the intervals of time in order to obtain a Strichartz inequality (for the frequency localized flow) on a fixed interval of time. Let us denote by $I=[0,1]$ and by $I_\delta=I+[-\delta, \delta]$, where $\delta$ is chosen like in (\ref{idPhi01}).

\begin{lem} Under the same notations as in Lemma \ref{1truncStr}, we have 
\begin{eqnarray}
\label{trunc2}
\nLdi{J_h^* \chi u_L}{I, L^6(M)} \leq & ch^{-\frac 34} \nLdi{J_h^* \chi u_L}{I_\delta, L^2(M)} +  \nonumber \\ & h^{\frac 34} \nLdi{\tilde{J}_h^* \tilde{\chi} u_L}{I_\delta, H^1(M)} + c\nLdi{\tilde{\chi} u_L}{I_\delta, H^1(M)}.
\end{eqnarray}
\end{lem}

\begin{proof} We sum the square of (\ref{trunc1}) over $t_0 \in \bar{J}$, where $\bar{J}$ was defined for the identity (\ref{idPhi01}). From (\ref{idPhi01}) and the definition of $\varphi$ we deduce that the reunion of intervals $I'_h(t_0)$, for $t_0\in \bar{J}$, recovers $I_\delta$ at most twice. Thus, $\sum_{t_0 \in \bar{J}} \nLdi{f}{I'_h}^2 \leq 2 \nLdi{f}{[-\delta, 1+\delta]}^2$. The result follows by merely observing that $\nLdi{J_h^* \chi u_L}{I, L^6(M)} \leq \nLdi{J_h^* \chi u_L}{I_\delta, L^6(M)}$ 
.
\end{proof}

From (\ref{trunc2}) we get the Strichartz inequality near the obstacle by means of Littlewood Paley summation.

\begin{propos} For every $\epsilon>0$ there exists $c_\epsilon>0$ such that for $(p,q)$ admissible in dimension $3$, 
\begin{equation}
\label{StrObpq}
\nX{\chi e^{it\tr_G} u_0}{L^p([0,1], W^{\frac {1}{2p}-\epsilon, q}(M))} \leq c_\epsilon \nX{u_0}{H^{\frac 1{p}} (M)}.
\end{equation}
\end{propos}

\begin{proof}
We apply a corollary of the Littlewood Paley theorem for $p,q\geq 2$ :
\begin{equation}
\nX{u}{L_T^p(W^{\sigma,q})} \leq c \nLL{S_0 u}{p}{q} + \left\{\sum_{j=0}^\infty 2^{2j\sigma} \nLL{\tr_j u}{p}{q}^2 \right\} ^{\frac 12}
\end{equation} 
Here we apply it for $(p,q)=(2,6)$ and $h=2^{-j}$, $\tr_j = J_{2^{-j}}^*$, $\sigma=\frac 14 -\epsilon >0$ and $u=\chi u_L$ on $I=[0,1]$. The left hand side term reads $\nX{\chi u_L}{L^2(I, W^{\sigma,6}(M))}$. Using (\ref{trunc2}), the parenthesis from the right hand side term is bounded by a sum $\sum_{j=0}^\infty$ of terms like 
$$ 2^{j(2-2\epsilon)} \nLdi{\tr_j \chi u_L}{I_\delta, L^2(M)}^2 + 2^{-j(1+2\epsilon)} \nLdi{\tr_j \chi u_L}{I_\delta, H^{1}(M)}^2 + 2^{-2j\epsilon} \nLdi{\tilde{\chi} u_L}{I_\delta, H^{1}(M)}^2.$$
Using the Plancherel theorem for the first two series and the geometric summation for the third (notice that $\sum_{j=0}^\infty 2^{-2j\epsilon} = c_\epsilon$), we obtain
$$\nX{\chi u_L}{L^2(I, W^{\sigma,6}(M))} \leq \nLdi{\chi u_L}{I_\delta, H^{1-\epsilon}(M)} + \nLdi{\chi u_L}{I_\delta, H^{\frac 12-\epsilon}(M)} + \nLdi{\tilde{\chi} u_L}{I_\delta, H^{1}(M)}.$$

We apply the smoothing effect (see Proposition 2.7 of \cite{BGTannIHP} and the translation onto the double). 
Thus, $$\nX{\chi u_L}{L^2(I, W^{\frac 14 - \epsilon,6}(M))} \leq c \nX{u_0}{H^{\frac 12}(M)}.$$ 
We want to perform a complex interpolation between the previous estimate and the conservation of the $L^2$ norm (we used also $0\leq \chi \leq 1$) : $$\nX{\chi u_L}{L^\infty(I, L^2(M))} \leq c \nLdi{u_0}{M}.$$
Using a weight of $\frac 2 p$, respectively $1-\frac{2}{p}$, we get an estimates of \Str type with loss of derivatives : 
$$\nX{\chi u_L}{L^p(I, W^{\frac 1{2p} - 2\epsilon,q}(M))} \leq c \nX{u_0}{H^{\frac 1{p}}(M)},$$ 
where $(p,q)$ satisfy $\frac{2}{p}+\frac{3}{q} = \frac{3}{2}$, i.e. they form an admissible couple in dimension $3$.
\end{proof}

\subsection{Proof of Proposition \ref{StrEstExtDom}}
\label{StrPr}

Combining estimates (\ref{StrObpq}) (\Str estimate near the boundary of $\Omega$) with (\ref{StrAway}) (\Str estimate away from the boundary)  for $s=\frac{1}{2p}-\epsilon$, we obtain, using that $\nX{v_0}{H^{\frac 1{2p}-\epsilon} (M)}\leq \nX{v_0}{H^{\frac 1{p}} (M)}$, 
$$\nX{e^{it\tr_G} v_0}{L^p([0,1], W^{\frac{1}{2p}-\epsilon, q}(M))} \leq c_\epsilon \nX{v_0}{H^{\frac 1{p}} (M)}.$$

Let $u_0\in H_D^{\frac {1}{2p} + \epsilon}(\Omega)$ and let $v_0\in H_{AS}^{\frac {1}{2p} + \epsilon}(\Omega)$ be such that ${v_0}_{|_\Omega}=u_0$. By uniqueness and stability at reflexion over the boundary of $\Omega$ of the linear flow (see Section \ref{Pre}), we have ${e^{it\tr_G}v_0}|_{\OO} = e^{it\tr_D} u_0$. Thus, $$\nX{e^{it\tr_G} v_0}{L^p([0,1], W^{s,q}(M))} \approx \nX{e^{it\tr_D} u_0}{L^p([0,1], W^{s,q}(\Omega))}$$ and $\nX{v_0}{H^{s} (M)} \approx \nX{u_0}{H^{s} (\Omega)}$. We obtain, 
$$\nX{e^{it\tr_D} u_0}{L^p([0,1], W^{\frac{1}{2p}-\epsilon, q}(\Omega))} \leq c_\epsilon \nX{u_0}{H^{\frac 1{p}} (\Omega)}.$$

We apply the ellipticity of the Laplacian $\tr_D$ to deduce a whole range of Strichartz inequalities : let $\tilde u_0 = (1-\tr_D)^{-\frac \sigma 2} u_0$, where $\sigma= \frac{1}{2p}-\epsilon - s$, $s\in [0,1]$. If $u_0 \in H^{s_0}(\OO)$ then $\tilde u_0 \in H^{s_0-\sigma}(\OO)$. We obtain the following inequality for $e^{it\tr_D} u_0$ :
\begin{equation}
\label{StrGM}
\nX{e^{it\tr_D} u_0}{L^p([0,1], W^{s,q}(\OO))} \leq c_\epsilon \nX{u_0}{H^{s+ \frac 1{2p} + \epsilon} (\OO)}.
\end{equation}

For $u_0\in H_N^{\frac{1}{2p}+\epsilon}(\Omega)$, we consider $v_0\in H_{S}^{\frac {1}{2p} + \epsilon}(M)$ be such that ${v_0}_{|_\Omega}=u_0$. We deduce as above the \Str inequality for the linear \Schr flow with Neumann Laplacian. 


\section{Global existence for NLS}
\label{NLSsect}
Having a \Str inequality we obtain classically a local existence theorem by Picard iteration scheme. These also enables propagation of the regularity of the initial data. Local existence in the energy space $H_0^1(\OO)$ combined with the conservation of the energy (and for defocusing nonlinearity of the $H_0^1(\OO)$ norm) enables us to conclude that the solution to is global in time.


\begin{proof}( of Theorem \ref{LeNLS})
Let us denote by $X_T=C([-T,T], H_0^1(\Omega)) \cap L^p([-T,T], L^\infty(\Omega))$ and, for a fix $u_0\in B\subset H_0^1(\Omega)$, by $\Phi : X_T \rightarrow X_T$ the functional
$$\Phi(u)(t) = e^{it\tr}u_0 -i \int_0^t e^{i(t-\tau)\tr}|u(\tau)|^2 u(\tau) \rm{d}\tau.$$
The space $X_T$ is a complete Banach space for the following norm
$$\nX{u}{X_T}={\max}_{|t|\leq T}\nX{u(t)}{H^1(\Omega)} + \nX{u}{L^p([-T,T],L^\infty(\Omega))}.$$ We prove that for a $T>0$ and $R>0$ small enough, $\Phi$ is a contraction from $B(0,R)\subset X_T$ into itself. We begin by estimating the $H^1$ norm of $\Phi(u)$ :
$$\nH{\Phi(u)(t)}{1} \leq \nH{u_0}{1} +  cT^{1-\frac {2}{p}} \nLL{u}{p}{\infty}^2 \nX{u}{L_T^\infty(H^1)} \leq \nH{u_0}{1} + c T^{1-\frac {2}{p}} \nX{u}{X_T}^{3}.$$
We have considered $2<p<3$.Thus, there exists $\epsilon>0$ such that $\epsilon<\frac{3}{2p}-\frac 12$. Therefore, by Sobolev imbedding theorem we have, for $(p,q)$ admissible in dimension $3$, that $W^{1-\frac{1}{2p}-\epsilon, q}(\Omega) \subset L^\infty(\Omega)$ : 
$$\nX{\Phi(u)}{L^p_T L^\infty}(\Omega) \leq c \nX{\Phi(u)}{L^p_T W^{1-\frac{1}{2p}-\epsilon, q}(\Omega)}.$$ 
Using the Strichartz estimate (\ref{StrExt}) and Minkowski inequality (like in the proof of (\ref{trunc1})), we have 
$$
\begin{array}{rcl}
\nX{\Phi(u)}{L^p(L^\infty)} & \leq & \nX{e^{it\tr}u_0}{L^p_T W^{1-\frac{1}{2p}-\epsilon, q}} + \nX{\int_0^t e^{i(t-\tau)\tr}|u|^2(\tau) u(\tau){\rm d}\tau} {L^p_T W^{1-\frac{1}{2p}-\epsilon, q}} \\
& \leq & c \nH{u_0}{1} + c \int_0^T \nX{|u|^{2}u(\tau)}{H^1(\Omega)}{\rm d}\tau.
\end{array}
$$
Using that $\nX{|u|^{2}(\tau)u(\tau)}{H^1(\Omega)} \leq c\nX{u(\tau)}{H^1} \nX{u(\tau)}{L^\infty}^2,$ we obtain
$$
\nX{\Phi(u)}{L^p(L^\infty)} \leq c \nH{u_0}{1} + cT^{1-\frac {2}{p}} \nX{u}{L^\infty(H^1)} \nX{u}{L^p(L^\infty)}^2 \leq c\nH{u_0}{1} + cT^{1-\frac {2}{p}} \nX{u}{X_T}^3.
$$
Thus, $\nX{\Phi(u)}{X_T} 
\leq c\nH{u_0}{1} + cT^{1-\frac {2}{p}} \nX{u}{X_T}^3.$

Consequently, there exist $T,R>0$, depending only on $B\subset H_0^1(\Omega)$ ($u_0\in B$), such that, for $u\in X_T$ with $\nX{u}{X_T}\leq R$, we have $\nX{\Phi(u)}{X_T}<R$.

As above, we prove that, for $u,v\in X_T$ such that $u(0)=u_0=v(0)$, $$\nX{\Phi(u) - \Phi(v)}{X_T} \leq c T^{1-\frac {2}{p}} (\nX{u}{X_T}^2 + \nX{v}{X_T}^2) \nX{u-v}{X_T}.$$

Choosing $T$ eventually smaller, we ensure that $\Phi$ is a contraction on the ball $B(0,R)\subset X_T$, $B(0,R)=\{u\in X_T, \ \nX{u}{X_T}<R \}$. Consequently, there exists a fix point of $\Phi$, which is therefore solution to (\ref{CubicNLS}).

For the Lipschitz property of the flow let us consider $u,v\in B(0,R)\subset X_T$ two solutions of (\ref{CubicNLS}) with initial data respectively $u_0, v_0 \in B$. As above, we have
$$\nX{u - v}{X_T} \leq c\nX{u_0 - v_0}{H^1} +c T^{1-\frac {2}{p}} (\nX{u}{X_T}^2 + \nX{v}{X_T}^2) \nX{u-v}{X_T}.$$ For $T,R>0$ chosen before we have $c T^{1-\frac {2}{p}} (\nX{u}{X_T}^2 + \nX{v}{X_T}^2)<1$ and therefore, $\exists \tilde{c}>0$ such that $\nX{u - v}{X_T} \leq \tilde{c} \nX{u_0 - v_0}{H^1}$. We conclude that the flow $u_0 \mapsto u$ is Lipschitz on $B\subset H_0^1$.

Let $\sigma\geq 1$ and suppose $u_0\in H^\sigma(\Omega) \cap H_0^1(\Omega)$. Let us estimate $\Phi(u)$ in $Y_T = C([-T,T], H^\sigma (\Omega)) \cap L^p([-T,T], L^\infty(\Omega))$ norm : $$\nX{u}{Y_T}={\max}_{|t|\leq T}\nX{u(t)}{H^\sigma(\Omega)} + \nX{u}{L^p([-T,T],L^\infty(\Omega))}.$$ As above, we obtain
$$\nX{\Phi(u)}{L_T^\infty H^\sigma} \leq c\nH{u_0}{\sigma} + cT^{1-\frac 2p} \nX{u}{X_T}^2 \nX{u} {L_T^\infty H^\sigma}.$$ We have chosen $T>0$ such that $c T^{1-\frac {2}{p}} (\nX{u}{X_T}^2 + \nX{v}{X_T}^2)<1$. Consequently, the $H^\sigma$ norm does not blow up for $|t|\leq T$ :
$$\nX{u}{L_T^\infty H^\sigma} \leq \tilde c \nH{u_0}{\sigma}.$$
Therefore we can conclude that regularity propagates.

The semilinear Schrödinger equation (\ref{CubicNLS}) has a Hamiltonian structure with gauge invariance and thus conservation laws hold for $H^2$ initial data. For $u_0\in H^1$ we deduce them by density :
the solution of (\ref{CubicNLS}) constructed above
satisfies, for $|t|\leq T$, to 
$$\left\{ \begin{array}{l}
\int |u(t)|^2 dx = \int |u_0|^2 dx, \\
\int |\nabla u(t)|^2 + \frac{1}{2}|u(t)|^{4} dx = \int |\nabla u_0|^2 + \frac{1}{2} |u_0|^{4} dx.
\end{array}
\right.$$

Moreover, note that $T>0$ depends only on $\nH{u_0}{1}$. Therefore, conservation of $H^1$ norm enables us to obtain, via a bootstrap argument, the global existence.
\end{proof}


\section{Global existence for Gross-Pitaevskii}
\label{GPsect}
The \GP equation (\ref{CubicGP}) is associated to the energy
\begin{equation}
\label{defEu}
{\mathcal E}(u)=\int_{\Omega} \frac{1}{2} |\nabla u|^2(x) +\frac{1}{4} \left( |u|^2(x)-1 \right)^2 dx.
\end{equation}
The main difference between the NLS (\ref{CubicNLS}) and the \GP equation (\ref{CubicGP}) is their energy space. For \GP it reads
$$E=\{u\in H_{loc}^1(\Omega),\ \nabla u \in L^2(\Omega),\ |u|^2-1 \in L^2(\Omega)\}.$$
Namely, the initial data in the energy space, $u_0\in E$, is not an $L^2(\Omega)$ function. Therefore we begin this section by describing the structure of $E$ and of the action of the linear \Schr group on $E$ by adapting the arguments of \cite{PG} to the boundary case. Then, we give the proof of the global existence theorem for the \GP equation (\ref{CubicGP}) by combining the latter structure with dispersive estimates derived in Section \ref{rappel} and \ref{near}.

\subsection{The energy space}
This section is inspired from \cite{PG}. In that paper, the Cauchy problem for \GP equation is studied in the whole Euclidean space $\RR^d$, for $d=2,3,4$. In the special case of $d=3$, $u_0\in E$ can be expressed in an explicit form as $u_0=c + v_0$, where $c\in \mathbb{C}$ and $v_0\in \dot{H}^1$. We show here that the same holds outside a non-trapping obstacle and give the outline of the proof. For more details we refer to \cite{PG}.

We denote by $C_0^\infty (\bar{\Omega})$ the restriction to $\bar{\Omega}$ of $C_0^\infty (\RR^3)$ and by $\pH$ the completion of  $C_0^\infty (\bar{\Omega})$ in the norm $\nLdi{\nabla \cdot}{\Omega}.$ We recall that $$\pH = \{u\in L^6(\Omega),\ \nabla u \in L^2(\Omega)\}.$$
Moreover, we have the following approximation property. 

Let $\chi \in C_0^\infty(\RR^3)$, $\chi=1$ on the ball of radius 1 $B(0,1)$ and $\chi=0$ outside $B(0,2)$. We define $\chi_R(x)=\chi(\frac{x}{R})$. For $v\in \pH$ we have $\chi_R v\in H^1(\OO)$ and 
\begin{equation}
\label{aproxPH}
\chi_R v \stackrel{R\rightarrow \infty}{\longrightarrow} {\rm v\ in\ the\ \nLdi{\nabla \cdot}{\Omega}\ norm.}
\end{equation}
%
%
%

We prove the main result of this section.

\begin{propos}
\label{EGP}
The energy space $E$ has the following structure
$$E=\{c+v,\ c\in \mathbb{C},\ |c|=1,\ v\in \pH,\ |v|^2 + 2Re(c^{-1}v) \in L^2(\Omega) \}.$$ The space $E$ is a complete metric space with the distance function
$$\delta_E(c+v, \tilde c + \tilde v) = |c-\tilde c| + \nLdi{\nabla v - \nabla \tilde v}{\Omega} + \nLdi{|v|^2 + 2Re(\bar{c} v) - |\tilde{v}|^2 - 2Re(\bar{\tilde{c}} \tilde{v})}{\Omega}.$$
\end{propos}
\begin{proof}
The embedding $"\supset"$ is obvious. For the converse we consider $R_0>0$ such that $\complement\Omega \subset B(R_0)$. For $u\in E$ we define, for every $\omega \in \mathbb{S}^2$ and $R>R_0$, $$U_R(\omega)=u(R\omega).$$ Just as in the proof of Lemma 7 of \cite{PG}, we show that $U_R$ 
converges to $U$ in $L^2(\mathbb{S}^2)$ norm and moreover $\nabla_\omega U = 0$. This enables us to conclude that $U$ is a constant $c(u)$. Since $|u|^2-1 \in L^2(\Omega)$, we conclude that $c(u)=1$. Let us proceed to the proof by noticing that 
\begin{equation}
\label{GradPol}
\int_{R_0}^\infty R^2 \nLdi{\partial_R U_R}{\Sph^2}^2 + \nLdi{\partial_\omega U_R}{\Sph^2}^2 dR \leq \nLdi{\nabla u}{\Omega}^2 < \infty .
\end{equation}
By Cauchy Schwarz, $\int_{R_0}^\infty \nLdi{\partial_R U_R}{\Sph^2} dR \leq \int_{R_0}^\infty R^2 \nLdi{\partial_R U_R}{\Sph^2}^2 dR$ and thus $\int_R^\infty \partial_\rho U_\rho d\rho$ satisfies the Cauchy criterion for convergence in $L^2(\Sph^2)$. We conclude the existence of a limit $U$ of $U_R$ in $L^2(\Sph^2)$. From (\ref{GradPol}) we deduce also that $\int_R^{R+1} \nLdi{\nabla_\omega U_\rho}{\Sph^2} d\rho$ goes to $0$ as $R\rightarrow \infty$. Since $\nabla_\omega U = \lim_{R\rightarrow \infty} \int_R^{R+1} \nabla_\omega U_\rho d\rho$ we conclude that $\nLdi{\nabla_\omega U}{\Sph^2} =0$. Thus, $U=c$, a constant of absolute value $1$.

Let us show that, if we denote by $v=u-c$, then $v\in \pH$. Notice that $\nabla v =\nabla u \in L^2(\Omega)$. Let $\chi \in C_0^\infty(\RR^3)$, $\chi=1$ on the ball of radius 1 $B(0,1)$ and $\chi=0$ outside $B(0,2)$. We define $\chi_R(x)=\chi(\frac{x}{R})$. We show that $v$ is the limit of $\chi_R v$ in the norm $\nLdi{\nabla \cdot}{\Omega}$. As $\chi_R v \in H^1(\OO)$, we obtain $v\in \dot{H}^1 (\Omega)$.

Notice that we have $v(R\omega) = -\int_R^\infty \partial_\rho U_\rho d\rho = -\int_R^\infty \omega \cdot (\nabla u)(\rho \omega) d\rho.$ By Cauchy Schwarz we obtain
$|v(R\omega)| \leq \frac{1}{\sqrt R} \left( \int_R^\infty \rho^2 |\nabla u|^2(\rho \omega) d\rho \right ) ^{\frac 12}.$ Consequently,
$$\int_{R'}^{2R'} \int_{\Sph^2} |v(R\omega)|^2 d\omega dR \leq \int_{R'}^{2R'} \frac {1}{R} \int_R^\infty \int_{\Sph^2} \rho^2 |\nabla u|^2(\rho \omega) d\omega d\rho dR.$$
Let us denote by $g(R)=\int_R^\infty \int_{\Sph^2} \rho^2 |\nabla u|^2(\rho \omega) d\omega d\rho$. The function $g$ is a decreasing function whose limit is $0$ at $\infty$. Then $\int_{R'}^{2R'} \frac {1}{R} g(R) dR < g(R') ln2$, which goes to $0$ as $R'$ goes to $\infty$. Consequently, $$\lim_{R'\rightarrow \infty} \int_{R'}^{2R'} \int_{\Sph^2} |v(R\omega)|^2 d\omega dR = 0.$$
This enables us to show that $\nLdi{\nabla \left( v-\chi_R v \right)}{\Sph^2} \rightarrow 0$ as $R\rightarrow \infty$. We have that $$\nabla \left( v-\chi_R v\right) = \frac{1}{R} (\nabla \chi)_R v + (1-\chi_R) v.$$ By writing $v$ in polar coordinates we obtain, for $R>R_0$,
$$\int_\Omega \frac 1{R^2} |(\nabla \chi)(\frac x R) v(x)|^2 dx \leq c\int_R^{2R} \int_{\Sph^2} |v(\rho \omega)|^2 d\omega d\rho \rightarrow 0$$ as $R\rightarrow \infty$. The other term also goes to 0 in $L^2(\Omega)$ norm as $R\rightarrow \infty$ : 
$$\nLdi{(1-\chi_R) \nabla v}{\Omega} \leq c \nLdi{\nabla v}{|x|>R} \rightarrow 0.$$
This concludes the proof of $v=u-c \in \pH$ and thus of the embedding $"\subset"$. The completeness of the metric space $E$ is an easy consequence of its structure.
\end{proof}

We end this section by showing that $E + H^1(\Omega) \subset E$ (see also Lemma 2 of \cite{PG}).

\begin{lem}
\label{Stab}
Let $u\in E$ and $w\in H^1(\Omega)$. Then $u+w \in E$ and
\begin{equation}
\label{EStabH1}
\nLdO{|u+w|^2 -1} \leq (\sqrt{{\mathcal E}(u)} + \nHO{w}{1})(1+\nHO{w}{1}).
\end{equation}
Moreover, for $\tilde u \in E$ and $\tilde w \in H^1(\Omega)$, we have
\begin{eqnarray}
\label{EDifH}
\delta_E(u+w, \tilde u + \tilde w) \leq (1+\nH{w}{1} + \nH{\tilde w}{1})\delta_E(u,\tilde u) + \nonumber \\ 
(1+ \sqrt{{\mathcal E}(u)} + \sqrt{{\mathcal E}(\tilde u)} + \nH{w}{1} + \nH{w}{1})\nH{w-\tilde w}{1}.
\end{eqnarray}
\end{lem}
\begin{proof}
From Proposition \ref{EGP} we know that $u=c+v$, $c\in \CC$, $|c|=1$ and $v\in \pH$. Then $u+w = c + (v+w)$ and $v+w \in \pH+ H^1(\Omega) \subset \pH$. We have to show that $v+w \in F_c$ or equivalent, that $|u+w|^2 - 1 \in L^2(\OO)$. We have $$|u+w|^2-1 = |v|^2 + 2Re(c^{-1}v) + |w|^2 + 2 Re(c^{-1} w) + 2 Re(\bar v w).$$
From Proposition \ref{EGP}  we have $|v|^2 + 2Re(c^{-1}v) \in L^2(\OO)$ and from (\ref{defEu}) $\nLdO{|v|^2 + 2Re(c^{-1}v)} \leq \sqrt{{\mathcal E}(u)}$. From $w\in H^1(\OO) \subset L^2(\OO)\cap L^6(\OO)$ we deduce $\nLdO{|w|^2}\leq c \nHO{w}{1}^2$, $\nLdO{2Re(\bar v w)} \leq c \nLi{v}{6}{\OO} \nHO{w}{1}$ and $\nLdO{2Re(\bar c w)} \leq c \nHO{w}{1}$. Estimate (\ref{EStabH1}) follows. For (\ref{EDifH}) we proceed similarly.
\end{proof}

\subsection{The action of $S(t)=e^{it\tr_N}$ on E}
\label{actStE}
This section is devoting to defining the action of the group $S(t)=e^{it\tr_N}$ on the energy space $E$. In view of the Neumann condition, $S(t)$ leaves constants invariant. We have to justify that $S(t)$ acts on $\pH$. We begin by recalling some functional calculus facts (e.g. \cite{RS}).

The domain of $-\tr_N$ in $L^2(\OO)$ is $H_N^2(\OO)=H^2(\OO)\cap\{\frac{\partial v}{\partial \nu} =0 \}$. For $v\in H_N^2(\OO)$ we have $\nLdO{\sqrt{-\tr_N}v} = \nLdO{\nabla v}$. Indeed, 
$$\nLdO{\sqrt{-\tr_N}v}^2 = (\sqrt{-\tr_N}v,\sqrt{-\tr_N}v)_{L^2} = (v, -\tr_N v)_{L^2} = \nLdO{\nabla v}^2.$$

The domain of $\sqrt{-\tr_N}$ in $L^2(\OO)$ is $H^1(\OO)$. For $u\in H^1(\OO)$ we also have the identity $\nLdO{\sqrt{-\tr_N}u} = \nLdO{\nabla u}$. Indeed, let $v\in H_N^2(\OO)$. Then
$$(\sqrt{-\tr_N}u,\sqrt{-\tr_N}v)_{L^2} = (u, -\tr_N v)_{L^2} = (\nabla u, \nabla v)_{L^2}.$$ From $\nLdO{\sqrt{-\tr_N}v} = \nLdO{\nabla v}$ for $v\in H_N^2(\OO)$ we deduce the same identity for $u\in H^1(\OO)$.

\begin{lem}
\label{SqN}
Using the notations of (\ref{aproxPH}), for $v\in \pH$ the limit $$lim_{R\rightarrow \infty} \sqrt{-\tr_N} (\chi_R v)$$ exists in the $L^2(\OO)$ norm and we denote it by $\sqrt{-\tr_N} v$. Moreover,
$$\nLdO{\sqrt{-\tr_N} v } = \nLdO{\nabla v}.$$
\end{lem}
\begin{proof}


From (\ref{aproxPH}) we have that $\left ( \nabla (\chi_R v) \right)_{R}$ is a Cauchy sequence in the $L^2(\OO)$ norm. As $\chi_R v\in H^1(\OO)$, the identity $\nLdO{\sqrt{-\tr_N}(\chi_R v)} = \nLdO{\nabla (\chi_R v)}$ holds. Therefore, $\left( \sqrt{-\tr_N}(\chi_R v) \right)_{R}$ is also a Cauchy sequence in the $L^2(\OO)$ norm. Denoting by $\sqrt{-\tr_N} v$ its limit, we obtain
$$\nLdO{\sqrt{-\tr_N} v } = \nLdO{\nabla v}.$$
\end{proof}
\begin{rem}
\label{funct}
Using the previous lemmas we can define a functional calculus $\varphi(\sqrt{-\tr_N})$ on $\pH$ for functions $\varphi : [0,\infty) \rightarrow \CC$ such that $\lambda\mapsto \frac{\varphi(\lambda)}{\lambda}$ is continuous and bounded for $\lambda \in [0,\infty)$. We denote by $$\varphi(\sqrt{-\tr_N}) v = \frac{\varphi(\sqrt{-\tr_N})}{\sqrt{-\tr_N}} \sqrt{-\tr_N}v$$ and this is well defined for $v\in \pH$ as $\sqrt{-\tr_N} v \in L^2(\OO)$. An equivalent definition is : $\varphi(\sqrt{-\tr_N}) v$ is the limit, in $L^2(\OO)$ norm, of $\varphi(\sqrt{-\tr_N}) (\chi_R v)$.
\end{rem}

An important consequence of the previous remark is the definition of $S(t)=e^{it\tr_N}$ on $\pH$. Let $v\in \pH$. We have $S(t)v = v + (e^{it\tr_N} - 1)v$ and each term of the sum is well defined.

\begin{lem}
\label{StpH}
For all $t\in \RR$ we have $S(t) : \pH \rightarrow \pH$ and moreover, for $v\in \pH$, we have
\begin{equation}
\label{SpHH}
\nHO{S(t)v -v}{1} \leq c(1+|t|^\frac 12) \nLdO{\nabla v}.
\end{equation}
\end{lem}
\begin{proof}
By functional calculus we have that $\frac{e^{it\tr_N}-1}{\sqrt{-\tr_N}} = \varphi(-\tr_N)$ acts on $L^2(\Omega)$ with a norm $\nX{\varphi(-\tr_N)}{L^2\rightarrow L^2} \leq sup_{\lambda \in \sigma(-\tr_N)}|\varphi(\lambda)|$. Here $\varphi(\lambda) = \frac{e^{it\lambda - 1}}{\sqrt \lambda}$, for $\lambda>0$. We have $\nL{\varphi}{\infty}\leq c \min(|t|\sqrt \lambda, \sqrt{\lambda ^{-1}})$. Optimising on $\lambda$ we obtain $\nL{\varphi}{\infty}\leq c|t|^\frac 12$ and thus 
$$\nLdi{\frac{e^{it\tr_N}-1}{\sqrt{-\tr_N}} \sqrt{-\tr_N} v}{\Omega} \leq c|t|^\frac 12 \nLdi{\nabla v}{\Omega}.$$
We have also $\nLdO{\sqrt{-\tr_N}(e^{it\tr_N}-1)v} \leq c \nLdO{\sqrt{-\tr_N} v} \leq c\nLdO{\nabla v}$.
Thus, $S(t)v =v + (S(t)-1)v \in \pH + H^1(\Omega) \subset \pH$.
\end{proof}

From the previous lemmas we shall deduce that $E$ is stable under the action of $S(t)$, for all $t\in \RR$.

\begin{propos}
For every $t\in \RR$ we have $S(t)E \subset E$. Moreover, for every $R>0$, for every $T>0$, there exists $C>0$ such that, for $u_0,\tilde u_0 \in E$ with ${\mathcal E}(u_0), {\mathcal E}(\tilde u_0) \leq R$, the following holds :
\begin{equation}
\label{StE}
sup_{|t|\leq T}\delta_E(S(t)u_0, S(t)\tilde u_0) \leq C \delta_E(u_0,\tilde u_0).
\end{equation}
\end{propos}
\begin{proof}
We write $S(t)u_0 = u_0 + (S(t)-1)u_0$. Writing $u_0=c_0+v_0$, with $v_0\in \pH$, we have that $S(t)u_0-u_0 = S(t)v_0 - v_0$.  From (\ref{SpHH}) we deduce $(S(t)-1)u_0 \in H^1(\OO)$. From Lemma \ref{Stab} we have $S(t)u_0 = u_0 + (S(t)-1)u_0 \in E$. Estimate (\ref{StE}) follows from (\ref{EDifH}), which reads in this setting : 
$$\delta_E(S(t)u_0, S(t)\tilde u_0) \leq c(1+|t|^\frac 12) (1+\sqrt{{\mathcal E}(u_0)} + \sqrt{{\mathcal E}(\tilde u_0)}) \delta_E(u_0, \tilde u_0).$$
\end{proof}

\subsection{\Str inequality and energy space}
As we mentioned in the beginning of Section \ref{GPsect}, one of the main differences between NLS and \GP is that the initial data is not in $L^2(\OO)$ for \GP. Therefore, it is not obvious to guess what the \Str inequality gives for $S(t)u_0$, when $u_0\in E$. This is the purpose of this section. We denote by $u_L(t)=S(t)u_0$, for all $t\in \RR$. We show in this section that for $u_0\in E$ and $2<p<3$ we have $u_L \in L^p([-T,T], L^\infty(\OO))$, for some $T>0$. We decompose $u_L$ in its high and low frequency parts and we treat them separately.

Let $\varphi_1\in C_0^\infty(\RR)$ such that $\varphi_1(s)=1$ pour $|s| \leq 1$ and $\varphi_1(s)=0$ pour $|s|\geq 2$. Let $\varphi_2\in C^\infty(\RR)$ such that $\varphi_1+\varphi_2 =1$. Let $u_0\in E$, $u_0= c_0 + v_0$, with $c_0\in \CC$, $|c_0|=1$ and $v_0\in \pH$. 

We denote by $v_{20} = \varphi_2(\sqrt{-\tr_N})v_0$. From Remark \ref{funct} and Lemma \ref{SqN} we deduce the following properties of $v_{20}$.

\begin{lem}
\label{v20}
Under the previous notations, we have $v_{20}\in H^1(\OO)$ and $$\nHO{v_{20}}{1} \leq c \nLdO{\nabla v_0}.$$
\end{lem}

In view of Lemma \ref{v20} we can apply the \Str inequality (\ref{StrExt}) (in Neumann setting) to $S(t)v_{20}$. 

\begin{lem}
Let $v_2(t)=S(t)v_{20}.$ For $T>0$ and $2<p<3$, the following holds : $v_2 \in L^p([-T,T], L^\infty(\OO))\cap L^\infty([-T,T], H^1(\OO))$ and $$\nX{v_2}{L_T^p(L^\infty)} + \nX{v_2}{L_T^\infty(H^1)} \leq C \nLdO{\nabla v_0}.$$
\end{lem}
\begin{proof}
From Lemma \ref{v20} we have $v_{20} \in H^1(\OO)$. Let $(p,q)$ be an admissible couple in dimension $3$ and $\epsilon>0$. From the \Str inequality (\ref{StrExt}) we deduce
$$\nX{v_2}{L_T^pW^{1-\frac{1}{2p}-\epsilon, q}(\OO)} \leq c \nHO{v_{20}}{1}.$$ For $2<p<3$ there exists $\epsilon>0$ such that $W^{1-\frac{1}{2p}-\epsilon, q}(\OO) \subset L^\infty(\OO)$ (see the proof of \ref{LeNLS}). Thus, $\nX{v_2}{L_T^pL^\infty(\OO)} \leq c \nLdO{\nabla v_0}.$ The estimate on $\nX{v_2}{L_T^\infty(H^1)}$ follows from the conservation of the $H^1$ norm by the linear \Schr flow $e^{it\tr_N}$.
\end{proof}

We denote by $v_{10}=v_0 - v_{20} = \varphi_1(\sqrt{-\tr_N})v_0$ and by $v_1(t) = S(t)v_{10}$.

\begin{lem}
\label{v10}
For $T>0$, there exists $C>0$ such that we have $v_1\in L^\infty([-T,T]\times \OO)$ satisfying $$\nX{v_1}{L_{t,x}^\infty} \leq C \nL{\nabla v_0}{2}.$$
\end{lem}
\begin{proof}In this proof we look at $v_1$ separately near the obstacle and away from the obstacle. The reason is that $v_{1}$ is only an $\pH$ function. Indeed, $\varphi_1(\sqrt{-\tr_N}) : L^6(\OO) \rightarrow L^6(\OO)$ and $\varphi_1(\sqrt{-\tr_N}) : L^2(\OO) \rightarrow L^2(\OO)$. As $S(t) : \pH \rightarrow \pH$ by Lemma \ref{StpH}, we obtain $v_{1}\in \pH$. 

We consider $\chi\in C_0^\infty(\RR^3)$ such that $\chi=1$ near $\Theta = \complement\OO$. Then $\chi v_1 \in L^\infty([-T,T],L^2(\OO))$ :  $$\nX{\chi v_1(t)}{L_T^\infty(L^2(\OO))} \leq \nLi{\chi}{3}{\OO} \nX{v_1}{L_T^\infty(L^6(\OO))} \leq C \nLdO{\nabla v_0}.$$
Similarly, we obtain $\tr(\chi v_1) = (\tr \chi) v_1 + 2 \nabla \chi \cdot \nabla v_1 + \chi (\tr v_1) \in L_T^\infty(L^2(\Omega))$. Moreover, $\frac{\partial}{\partial \nu}(\chi v_1)|_{\partial\OO} = \frac{\partial v_1}{\partial \nu}|_{\partial\OO}=0 $ as $\chi=1$ in the neighborhood of $\partial\OO$. Thus, $\chi v_1 \in L_T^\infty H_N^2(\OO)$, where $H_N^2(\OO)$ is the domain of $-\tr_N$ in $L^2(\OO)$. As $H_N^2(\OO) \subset L^\infty(\OO)$, we obtain $\chi v_1 \in L^\infty([-T,T]\times\OO)$.

We pass to the term $(1-\chi)v_1$. It can be seen as a function on $\RR^3$ in the $x$ variable extending it by $0$. Since $v_1 \in L_T^\infty L^6(\OO)$, we have $(1-\chi) v_1 \in L_T^\infty L^6(\RR^3)$. We show that $(1-\chi)v_1 \in L_T^\infty W^{2,6}(\RR^3)$. For that purpose, it suffices to show that $\tr((1-\chi)v_1) \in L_T^\infty(L^6(\RR^3))$. We have
\begin{equation}
\label{calcul}
\tr((1-\chi)v_1) = -(\tr \chi)v_1 -2\nabla \chi\cdot \nabla v_1 + (1-\chi)(\tr v_1).
\end{equation}
Clearly, the first and the last term of the right hand side expression are in $L_T^\infty(L^6(\RR^3))$. For $\nabla \chi\cdot \nabla v_1$ we need to do finer analysis. As $\nabla v_1 \in L_T^\infty L^2(\OO)$ we deduce $\nabla \chi\cdot \nabla v_1 \in L_T^\infty L^2(\RR^3)$. We show that $\nabla \chi\cdot \nabla v_1 \in L_T^\infty W^{2,2}(\RR^3).$ We compute 
$$\tr(\nabla \chi\cdot \nabla v_1) = (\tr \nabla \chi)\cdot \nabla v_1 + 2 (\nabla^2 \chi) \cdot (\nabla^2 v_1) + \nabla \chi\cdot (\tr \nabla v_1).$$ 
We have $(\tr \nabla \chi)\cdot \nabla v_1 \in L_T^\infty L^2(\RR^3)$ and $\nabla \chi\cdot (\tr \nabla v_1) \in L_T^\infty L^2(\RR^3)$. The middle term, $2 (\nabla^2 \chi) \cdot (\nabla^2 v_1)$ can be written as $P(x,D)(1-\tr)v$, with $P(x,D) = 2 (\nabla^2 \chi) \cdot (\nabla^2 (1-\tr)^{-1})$ an pseudo-differential operator of order $0$ with compact support. Its coefficients are independent of $t$. Consequently, $2 (\nabla^2 \chi) \cdot (\nabla^2 v_1) \in L_T^\infty L^6(\RR^3)$ and since this function is compactly supported in $x$, it belongs also to $L_T^\infty L^2(\RR^3).$ 

We obtain $\nabla \chi\cdot \nabla v_1 \in L_T^\infty W^{2,2}(\RR^3) \subset L_T^\infty L^6(\RR^3).$ Going back to (\ref{calcul}) we deduce $(1-\chi)v_1 \in L_T^\infty W^{2,6}(\RR^3) \subset L^\infty([-T,T] \times\RR^3)$. Taking the restriction to $\OO$ concludes the proof.
\end{proof}

From the previous lemmas, we deduce easily the following.

\begin{propos}
\label{propLL}
For $T>0$ and $2<p<3$, there exists $C>0$ such that, for $u_0\in E$ and $u_L(t)=e^{it\tr_N} u_0$, we have : 
$u_L \in L^p([-T,T], L^\infty(\OO))$ and 
\begin{equation}
\label{estuLLpL}
\nX{u_L}{L_T^p(L^\infty)} \leq 1 + C \nLdO{\nabla v_0}.
\end{equation}
Moreover, for $\tilde u_0 \in E$ and $\tilde u_L(t)=e^{it\tr_N} \tilde u_0$,
\begin{equation}
\label{estDifuLLpL}
\nX{u_L - \tilde u_L}{L_T^p(L^\infty)} \leq C \delta_E(u_0,\tilde u_0).
\end{equation}
\end{propos}
\begin{proof}
We write $u_L(t)= c_0 + e^{it\tr_N} v_0 = c_0 + v_1(t) + v_2(t)$. The conclusion follows from  $c_0 \in \CC$, $v_1 \in L_T^p(L^\infty)$, $v_2 \in L_T^\infty(L^\infty)$ and their respective estimates.
\end{proof}

We close this section by collecting estimates which will be useful in the sequel.
We consider $u_0, \tilde u_0\in E$, $u_L(t)=S(t)u_0$ and $\tilde u_L(t)=S(t) \tilde u_0$, $w,\tilde w \in X_T=C([-T,T], H_0^1(\Omega)) \cap L^p([-T,T], L^\infty(\Omega))$ with the associated norm $\nX{w}{X_T}={\max}_{|t|\leq T}\nX{w(t)}{H^1(\Omega)} + \nX{w}{L^p([-T,T],L^\infty(\Omega))}.$ Let $u=u_L+w$ and $\tilde u = \tilde u_L + \tilde w$. We denote by $$\gamma = \delta_E(u_0, \tilde u_0) + \nX{w - \tilde w}{X_T}.$$

\noindent As a corollary of Lemma \ref{StpH} and \ref{Stab} we have
\begin{equation}
\label{1} \nX{|u_L|^2 - 1}{L_T^\infty L^2(\OO)} \leq c(\sqrt{{\mathcal E}(u_0)} + {\mathcal E}(u_0))
\end{equation}
\begin{equation}
\label{2} 
\nX{|u|^2 - 1}{L_T^\infty L^2(\OO)} \leq (1 + {\mathcal E}(u_0))\nX{w}{X_T} + \nX{w}{X_T}^2.
\end{equation}

\noindent As a corollary of Proposition \ref{propLL} we have
\begin{equation}
\label{3} \nX{u}{L_T^p L^\infty(\OO)} \leq c(1 + \sqrt{{\mathcal E}(u_0)} + \nX{w}{X_T})
\end{equation}
\begin{equation}
\label{4} 
\nX{u - \tilde u}{L_T^p L^\infty(\OO)} \leq \gamma
\end{equation}
From (\ref{3}), (\ref{4}) and (\ref{StE}) we deduce
\begin{equation}
\label{5} 
\nX{|u|^2 - |\tilde u|^2 }{L_T^\infty L^2} \leq \gamma (1+ \sqrt{{\mathcal E}(u_0)} + \sqrt{{\mathcal E}(\tilde u_0)}+ \nX{w}{X_T} + \nX{\tilde w}{X_T}) 
\end{equation}
Moreover,
\begin{equation}
\label{9} 
\nX{|u|^2 - 1 }{L_T^{\frac p2} L^\infty} \leq 1+ {\mathcal E}(u_0) + \nX{w}{X_T}^2
\end{equation}
\begin{equation}
\label{8} 
\nX{|u|^2 - |\tilde u|^2 }{L_T^{\frac p2} L^\infty} \leq \gamma (1+ \sqrt{{\mathcal E}(u_0)} + \sqrt{{\mathcal E}(\tilde u_0)} + \nX{w}{X_T} + \nX{\tilde w}{X_T})
\end{equation}

\noindent By simple computations we obtain
\begin{equation}
\label{7} 
\nX{\nabla u}{L_T^\infty L^2(\OO)} \leq \sqrt{{\mathcal E}(u_0)} + \nX{w}{X_T}
\end{equation}
\begin{equation}
\label{6} 
\nX{\nabla u - \nabla \tilde u}{L_T^\infty L^2(\OO)} \leq \gamma
\end{equation}

The estimates (\ref{1}) to (\ref{6}) follow from simple computations, decomposing $u=u_L+w$ and applying H\"older and Sobolev inequalities combined with the estimates cited.

\subsection{Proof of Theorem \ref{GeGP}}
\label{sectProofGP}

Let $u_0\in E$. In Section \ref{actStE} we presented the action of $S(t)=e^{it\tr_N}$ on $E$. We recall the notation $u_L(t)=S(t)u_0$. We call the solution of (\ref{CubicGP}) the solution to the Duhamel associated formula :
\begin{equation}
\label{GPDuhamel}
u(t)=u_L(t) - i \int_0^t e^{i(t-\tau)\tr_N} F(u)(\tau) d\tau,
\end{equation}
where $F(u)=(|u|^2 - 1)u$. We denote by $w=u-u_L$ and by $\Phi$ the functional
\begin{equation}
\label{defPhiGP}
\Phi(w)= - i \int_0^t e^{i(t-\tau)\tr_N} F(u_L + w)(\tau) d\tau.
\end{equation}

We show the local existence of $u$ that satisfies (\ref{GPDuhamel}) by showing that $\Phi$ has a fixed point $\Phi(w)=w$. For that purpose we define, for $T>0$ and $2<p<3$, $X_T=C([-T,T], H_0^1(\Omega)) \cap L^p([-T,T], L^\infty(\Omega))$. The space $X_T$ is a complete Banach space for the following norm
$$\nX{w}{X_T}={\max}_{|t|\leq T}\nX{w(t)}{H^1(\Omega)} + \nX{w}{L^p([-T,T],L^\infty(\Omega))}.$$ We prove that, for a $T>0$ and $R>0$ small enough, $\Phi$ is a contraction from $B(0,R)\subset X_T$ into itself.

\begin{lem}
\label{PhiXT}
Using the previous notations we have, for $w\in X_T$, that $$\nX{\Phi(w)}{X_T} \leq c\nX{F(u_L + w)}{L_T^1 H^1(\OO)}.$$
\end{lem}
\begin{proof}
From (\ref{GPDuhamel}) we deduce, by Minkowski inequality, that $$\nX{\Phi(w)}{L_T^\infty L^2(\OO)} \leq c \nX{F(u_L+w)}{L_T^1 L^2(\OO)}.$$ As $\nabla(\Phi(w))= - i \int_0^t e^{i(t-\tau)\tr_N} \nabla (F(u_L + w)) (\tau) d\tau$ we have also $$\nX{\nabla(\Phi(w))}{L_T^\infty L^2(\OO)} \leq c \nX{\nabla(F(u_L+w))}{L_T^1 L^2(\OO)}.$$ We have considered $2<p<3$. Thus, there exists $\epsilon>0$ such that, for $(p,q)$ an admissible couple in dimension $3$, $W^{1-\frac{1}{2p}-\epsilon, q}(\OO) \subset L^\infty(\OO)$. From the Strichartz inequality (\ref{StrExt}) we obtain :
$$\nX{\Phi(w)}{L_T^p L^\infty(\OO)} \leq \nX{\Phi(w)}{L_T^p W^{1-\frac{1}{2p}-\epsilon, q}(\OO)} \leq c \nX{F(u_L+w)}{L_T^1 H^1(\OO)}.$$
\end{proof}

We have to estimate $F(u)$ in $L_T^1 H^1(\OO)$ for $u=u_L + w$, $w\in X_T$. For the fixed point method we also need to estimate $\nX{F(u_L + w) - F(\tilde u_L + \tilde w)}{L_T^1 H^1(\OO)}$.

\begin{propos}
\label{estF}
Under the conditions of Section \ref{sectProofGP} we have
\begin{eqnarray*}
& \nX{F(u)}{L_T^1L^2} \leq c T^{1-\frac 1p} (1+ {\mathcal E}(u_0) + \nX{w}{X_T})^2 \nX{w}{X_T} \\
& \nX{\nabla (F(u))}{L_T^1L^2} \leq c T^{1-\frac 2p} (1+ \sqrt{{\mathcal E}(u_0)} + \nX{w}{X_T})^3  
\end{eqnarray*}
and
\begin{eqnarray*}
& \nX{F(u) - F(\tilde u)}{L_T^1L^2} \leq c T^{1-\frac 1p} \gamma (1+ {\mathcal E}(u_0)+ {\mathcal E}(\tilde u_0) + \nX{w}{X_T}+ \nX{\tilde w}{X_T})^2 \\
& \nX{\nabla (F(u) - F(\tilde u))}{L_T^1L^2} \leq c T^{1-\frac 2p} \gamma (1+ \sqrt{{\mathcal E}(u_0)} + \sqrt{{\mathcal E}(\tilde u_0)}+ \nX{w}{X_T} + \nX{\tilde w}{X_T})^2
\end{eqnarray*}
where we have denoted by $\gamma = \delta_E(u_0, \tilde u_0) + \nX{w - \tilde w}{X_T}$.
\end{propos}
Notice that, if $u_0 = \tilde u_0$, then we have $\gamma = \nX{w - \tilde w}{X_T}$.
\begin{proof}
The conclusions follow from estimates (\ref{1}) to (\ref{6}). Let us explain one of the conclusions, for example the estimate on $F(u_L+w) - F(\tilde u_L + \tilde w)$. We have
$$F(u) - F(\tilde u) = (|u|^2- |\tilde u|^2)u + (u-\tilde u)(|u|^2-1).$$ We apply the H\"older inequality combined with (\ref{5}) and (\ref{3}) for the first term and (\ref{4}) and (\ref{2}) for the second one. We bound thus $\nX{F(u_L+w) - F(\tilde u_L + \tilde w)}{L_T^{p} L^2}$. By H\"older inequality we obtain the positive power of $T$ : $$\nX{F(u) - F(\tilde u)}{L_T^1L^2} \leq cT^{1-\frac 1p}\nX{F(u) - F(\tilde u)}{L_T^p L^2}.$$ The other estimates follow similarly.
\end{proof}

Combining the estimates on the nonlinear term from Proposition \ref{estF} with Lemma \ref{PhiXT} we obtain the following
\begin{cor} Under the conditions on Lemma \ref{PhiXT} we have
\begin{eqnarray}
\label{estPhi}
& \nX{\Phi(w)}{X_T} \leq c T^{1-\frac 2p} (1+ {\mathcal E}(u_0) + \nX{w}{X_T})^3\\
\label{estDifPhi}
& \nX{\Phi(w) - \Phi(\tilde w)}{X_T} \leq c T^{1-\frac 2p} \gamma (1+ {\mathcal E}(u_0)+ {\mathcal E}(\tilde u_0) + \nX{w}{X_T}+ \nX{\tilde w}{X_T})^2,
\end{eqnarray}
 where we denoted by $\gamma = \delta_E(u_0, \tilde u_0) + \nX{w - \tilde w}{X_T}$.
 \end{cor}

As a consequence, we can prove the global wellposedness result from Theorem \ref{GeGP} on \GP equation (\ref{CubicGP}).

\begin{proof}We fix $u_0\in B \subset E$. 
From estimate (\ref{estPhi}) we deduce that there exist $T,R>0$, depending only on $B\subset E$ ($u_0\in B$), such that, for $w\in X_T$ with $\nX{w}{X_T}\leq R$, we have $\nX{\Phi(w)}{X_T}<R$.

For $\tilde u_0 = u_0$ estimate (\ref{estDifPhi}) reads  $$\nX{\Phi(w) - \Phi(\tilde w)}{X_T} \leq c T^{1-\frac 2p}(1+ \nX{w}{X_T}+ \nX{\tilde w}{X_T})^2 \nX{w - \tilde w}{X_T}.$$ As $2<p$, choosing $T$ eventually smaller ensures that $\Phi$ is a contraction on the ball $B(0,R)\subset X_T$, $B(0,R)=\{w\in X_T, \ \nX{w}{X_T}<R \}$. Consequently, there exists a fixed point of $\Phi$ in $B(0,R)$, which is therefore solution to (\ref{CubicGP}).

For the Lipschitz property of the flow let us consider $u, \tilde u\in B(0,R)\subset X_T$ two solutions of $\Phi(u-u_L) = u-u_L$, therefore of (\ref{CubicGP}), with initial data respectively $u_0, \tilde u_0 \in B$. 

From (\ref{estDifPhi}) we have, for $w=u-u_L$ and $\tilde w= \tilde u- \tilde u_L$, 
$$\nX{w - \tilde w}{X_T} \leq c T^{1-\frac 2p} (\delta_E(u_0, \tilde u_0) + \nX{w - \tilde w}{X_T}) (1+ {\mathcal E}(u_0)+ {\mathcal E}(\tilde u_0) + \nX{w}{X_T}+ \nX{\tilde w}{X_T})^2,$$

For $T,R>0$ chosen before we have $c T^{1-\frac {2}{p}} (1+ {\mathcal E}(u_0)+ {\mathcal E}(\tilde u_0) + \nX{w}{X_T}+ \nX{\tilde w}{X_T})^2<1$ and therefore, $\exists \tilde{c}>0$ such that $$\nH{w-\tilde w }{1} \leq \nX{w - \tilde w}{X_T} \leq \tilde{c} \delta_E(u_0, \tilde u_0).$$ 
From (\ref{EDifH}) we have $\delta_E(u(t),\tilde u(t)) \leq C(R,B) (\delta_E(u_0, \tilde u_0) + \nX{w -\tilde w}{L_T^\infty H^1})$. Consequently, there exists $C>0$ such that $\delta_E(u(t), \tilde u(t)) \leq c\delta_E(u_0, \tilde u_0)$, for all $|t|\leq T$. We conclude that the flow $u_0 \mapsto u(t)$ is Lipschitz on $B\subset E$.

The proof of the propagation of regularity from section 3.3 of \cite{PG} adapts to the framework of exterior domains using techniques similar to those of Section \ref{actStE}. Those techniques combined with the stability of $E$ by summation with an $H^1$ element (see Lemma \ref{Stab}) enables us to show that $u_0\in E$ can be approached, in $\delta_E$ distance, by $u_0^\epsilon \in E$ such that $\tr u_0^\epsilon \in L^2(\OO)$. As one can prove conservation of energy ${\mathcal E}$ for initial data $f\in E$ such that $\tr f \in L^2(\OO)$, from (\ref{StE}) we deduce that conservation of energy holds for $u_0\in E$ : ${\mathcal E}(u(t)) = {\mathcal E}(u_0)$.

Notice that $T$, the existence time for which we applied the fixed point method, depends on ${\mathcal E} (u_0)$ and on $R$. From the conservation of energy for the solutions of (\ref{CubicGP}) we have ${\mathcal E} (u_0) = {\mathcal E} (u(t))$ for all $|t|\leq T$. Consequently, we can apply a bootstrap argument and conclude to the extension globally in time of $u \in C(\RR, E)$, solution of (\ref{CubicGP}).

\end{proof}


\end{document}